\documentclass{amsart}
\usepackage{amssymb}
\usepackage{graphicx}
\usepackage{xcolor}
\usepackage{ulem}
\usepackage[colorlinks=true,citecolor=black,linkcolor=black,urlcolor=blue]{hyperref}

\numberwithin{equation}{section}

\newtheorem{theorem}{Theorem}

\theoremstyle{definition}

\theoremstyle{remark}

\newtheorem{example}[theorem]{Example}
\newcommand{\bExample}{\begin{example}}
\newcommand{\eExample}{\end{example}}

\newcommand{\benum}{\begin{enumerate}}
\newcommand{\eenum}{\end{enumerate}}
\newcommand{\bitemize}{\begin{itemize}}
\newcommand{\eitemize}{\end{itemize}}
\newcommand{\bEqn}{\begin{equation}}
\newcommand{\eEqn}{\end{equation}}

\newcommand{\rr}{\mathbb R}
\newcommand{\bd}{\partial}

\newcounter{arb}

\newcommand{\bibtitle}[1]{{#1}}

\title{A Gap in Stanfield's Proof of Sachs' Linear Linkless Embedding Conjecture}
\author{Ramin Naimi}
\address{Department of Mathematics,
Occidental College,
Los Angeles, CA 90041, USA.}

\date \today

\begin{document}{


\maketitle

The aim of this note is to describe what I believe is a serious gap 
in Stanfield's~\cite{LS} proof of Sachs' conjecture
that
every linklessly embeddable graph 
has a  linear linkless embedding in $\rr^3$.
In \cite{LS}, page~4594, lines 5-6, it says:
\begin{quote}
\textit{
$(*)$ 
Because $x$ is so close to the position of $v$, 
$\Gamma'$  is also disjoint from any edges incident
to $x$, since it was disjoint from the edges of $v$.
}
\end{quote}
Below, I give an example that
contradicts statement $(*)$ above.
Before doing so I will give a brief description of
the argument  in \cite{LS}
and what $(*)$ is referring to.

Let $G$ be a linklessly embeddable graph.
Let $\phi$ be a paneled embedding of $G$ in $\rr^3$
(this means every cycle in $\phi$ 
bounds a disk in whose interior is disjoint from $\phi$).
We wish to show by induction on the size of $G$
that $\phi$ is ambient isotopic to a linear embedding.
Let $xy$ be an edge of $\phi$.
Then, by the induction hypothesis,
there is an ambient isotopy $F: \rr^3 \times I \to \rr^3$ such that 
$F_0$ is the identity map on $\rr^3$ and
$F_1(\phi/xy)$ is a linear embedding.

Let $B_{xy}$ be a ball that is a closed $\epsilon$-neighborhood
of $xy \subset \phi$.
Let $v$ be a point in the interior of $xy$.
Let $D_{xy}$ be a disk in $B_{xy}$ such that 
$D_{xy} \cap \phi = D_{xy} \cap xy = \{v\}$ and
$D_{xy} \cap \bd B_{xy} = \bd D_{xy}$.
We contract $xy$ to $v$,
keeping $\rr^3 \setminus B_{xy}$ pointwise fixed.
Let $D_F := F_1(D_{xy})$.
(In \cite{LS}
this disk is denoted as just $D_{xy}$.)

Now we try to modify $F_1(\phi/xy)$ to construct a linkless linear embedding $\Psi$ of $G$.
We place the vertices $x$ and $y$ very close to, and on opposite sides of, $D_F$.
For each neighbor $w$ of $v$ in $G/xy$,
the straight edge $vw$ in $F_1(\phi/xy)$ 
is replaced by the straight edge(s) $xw$ or $yw$ or both
according to whether 
$w$ is adjacent to $x$, $y$, or both, in $G$.

To show the linear embedding $\Psi$ is paneled,
part of the argument in \cite{LS} goes as follows.
Suppose $C$ is a cycle in $\Psi$
that contains $x$ but not $y$.
Since $F_1(\phi/xy)$ is paneled, 
the cycle $C' := C /xy$ 
bounds a disk $\Gamma'$ whose interior is disjoint from $F_1(\phi/xy)$.
It is also shown in \cite{LS} 
that the interior of $\Gamma'$  is disjoint from $D_F$;
thus $\Gamma' \cap D_F = \{v\}$ $\subset \bd \Gamma'$.

\smallskip

Now here is the gap.
Statement~$(*)$ above claims the interior of $\Gamma'$ is disjoint from
edges incident to $x$.
But the following construction and embedding of a graph $G$
and disks $\Gamma'$ and $D_F$ contradict $(*)$.

Let $G$ be the graph shown in Figure~\ref{fig-1},
with $n$ to be defined later.
(This graph is planar and therefore has a linear linkless embedding.
We use $G$ only to illustrate
why $(*)$ is not a valid statement,
not as a counterexample to the claim
that $\psi$ is linkless,
which may or may not be true in general.)
In Figure~\ref{fig-2}, $S^2$ is a sphere in $\rr^3$ centered at $v$;
let $\alpha \subset S^2$ be the simple curve shown in the figure.
Let $\Delta$ be the disk obtained by coning $v$ onto $\alpha$ 
with straight line segments
(coning the center, $v$, of $S^2$ onto any simple curve in $S^2$ always yields a disk).
We embed the neighbors $a_1, b_1, c$ of $x$ in $S^2$ as shown in Figure~\ref{fig-2}.

{
\begin{figure}[h]
\includegraphics[width=100mm]{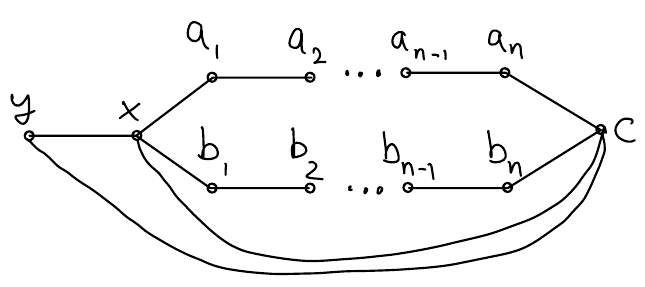}
\caption{
\label{fig-1}
}
\end{figure}

\begin{figure}[h]
\includegraphics[width=95mm]{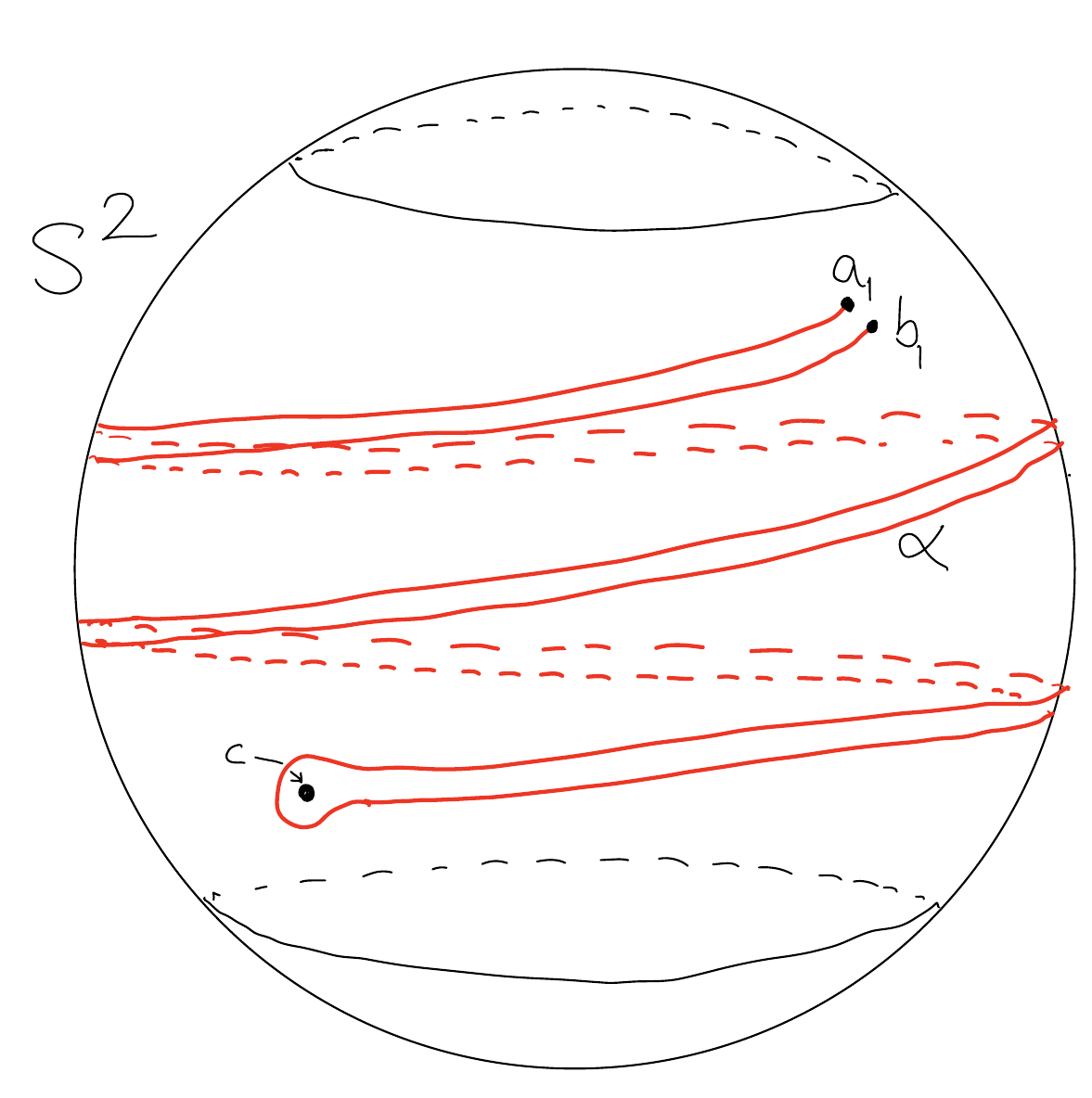}
\caption{
\label{fig-2}
}
\end{figure}

Then, regardless of how close $x$ is to $v$,
at least one of the line segments $a_1 x, b_1 x, cx$
will intersect the interior of $\Delta$
(because
if $x$ is placed anywhere off the center of $S^2$ 
such that $cx$ is disjoint from $\Delta$,
then, for every point $p \in S^2$ above the equator,
the line segment $px$ intersects the interior of $\Delta$).

Let
$\eta, \eta' \subset S^2$ be the arcs shown in Figure~\ref{fig-3}.
Let
$\alpha' $
$:= (\alpha  \setminus \eta) \cup \eta'$,
and $\Delta' := \Delta \cup \delta$,
where $\delta \subset S^2$ is the disk bounded by $\eta \cup \eta'$ as shown.
Then the interior of $\Delta'$ 
intersects at least one of the line segments $a_1 x, b_1 x, cx$ (as does $\Delta$).
Now we embed the vertices
$a_2, \cdots, a_n$
along the arc $a_1 c \subset \alpha'$, and
$b_2, \cdots, b_n$
along the arc $b_1 c \subset \alpha'$,
as shown in Figure~\ref{fig-3}.
By making $n$ large enough, 
we can make the line segments 
$a_i a_{i+1}$, $b_i b_{i+1}$, $a_n c$, and $b_n c$
short enough to make the linear path 
$a_1 \cdots a_n c b_n \cdots b_1$
arbitrarily close to $\alpha'$.
Let $\Gamma' \subset \Delta'$
be the disk bounded by the linear cycle
$v a_1 \cdots a_n c b_n \cdots b_1 v$.
Then, no matter where $x$ is,
$\Gamma'$ intersects at least one of $a_1 x, b_1 x, cx$.

\begin{figure}[h]
\includegraphics[width=95mm]{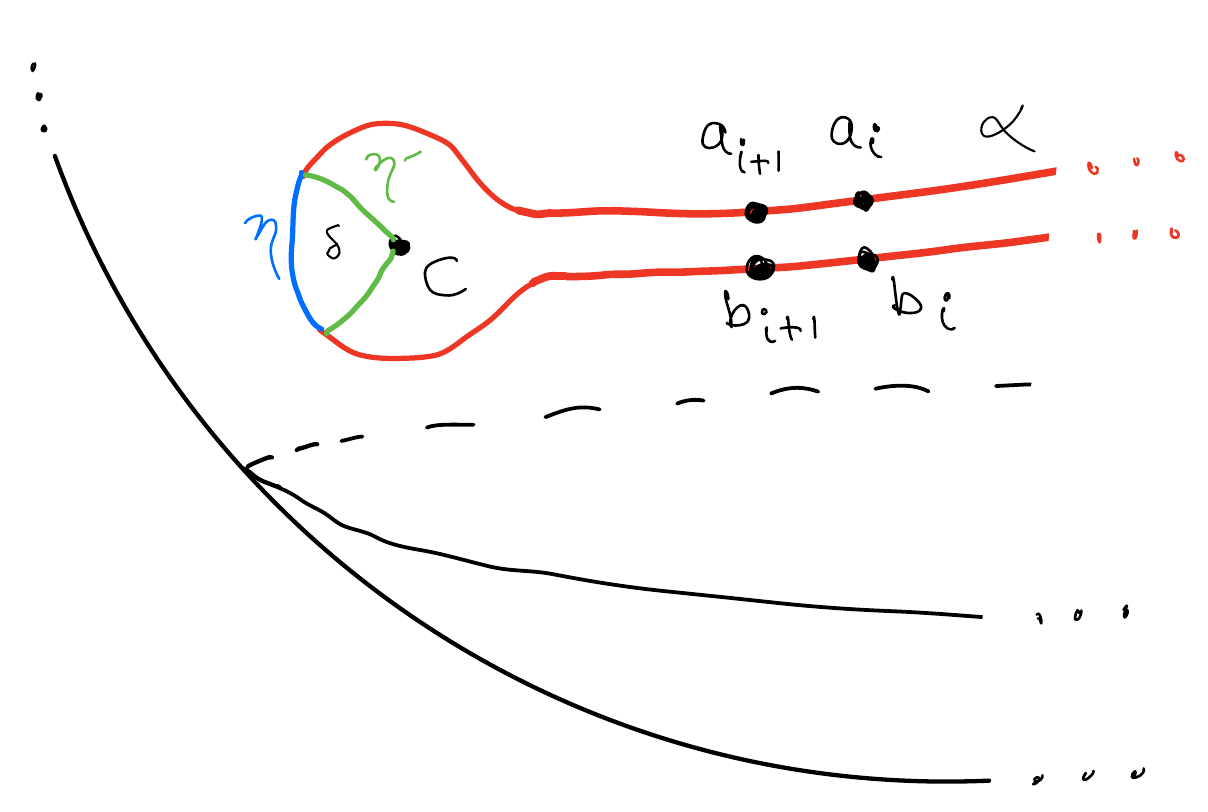}
\caption{
\label{fig-3}
}
\end{figure}

Lastly, we show the disk $D_F$ can be embedded in a manner
compatible with the argument given in \cite{LS}.
Let $\beta \subset S^2$ be the simple closed curve 
shown in Figure~\ref{fig-4}.
Let $D_F$ be the disk obtained by coning $v$
onto $\beta$.
Then $D_F \cap \Gamma' = \{v\}$
and $D_F$ is disjoint from the interior of $\Gamma'$
since $\beta$ is disjoint from $\alpha'$.
So this construction contradicts statement $(*)$.

Note that although we can assume $D_{xy}$ is a ``flat'' disk,
we cannot assume it is flat 
after the ambient isotopy $F$.
Perhaps it was such an assumption that 
led to~$(*)$.

\begin{figure}[h]
\includegraphics[width=95mm]{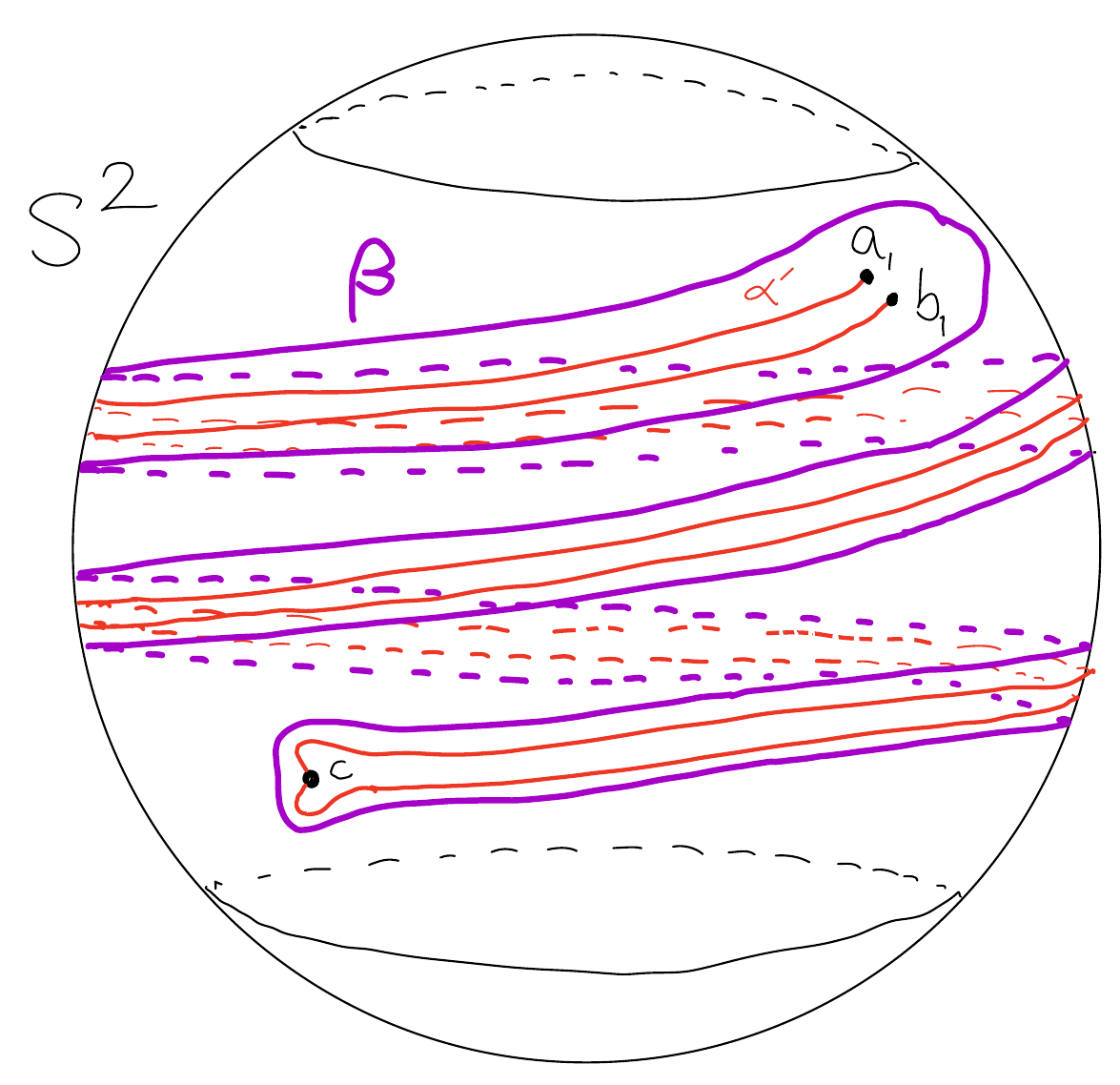}
\caption{
\label{fig-4}
}
\end{figure}

\end{document}